\documentclass[12pt]{amsart}
\usepackage{epsfig,comment}
\usepackage{mathtools}
\usepackage[subnum]{cases}
\usepackage{enumerate}
\usepackage{bbm,bm}
\usepackage{amsmath,amssymb,mathrsfs,amsthm}
\usepackage{color}
\usepackage[bookmarks=true,
bookmarksnumbered=true, breaklinks=true,
pdfstartview=FitH, hyperfigures=false,
plainpages=false, naturalnames=true,
colorlinks=true,pagebackref=true,
pdfpagelabels,colorlinks=true,
	linkcolor=black,
	citecolor=black,
	filecolor=black,
	urlcolor=black]{hyperref}
\usepackage{wasysym}
\usepackage[hang]{footmisc}
\usepackage{verbatim}

\usepackage{tikz}
\usetikzlibrary{calc,arrows.meta,decorations.pathreplacing}
\usetikzlibrary{intersections}

\newcommand{\E}{\mathbb{E}}

\newcommand{\N}{\mathbb{N}}

\newcommand{\R}{\mathbb{R}}

\newcommand{\NN}{N^{-\frac{2}{d-1}}}
\newcommand{\NNN}{N^{\frac{2}{d-1}}}

\DeclareMathOperator{\vol}{vol}
\DeclareMathOperator{\Gr}{Gr}
\DeclareMathOperator{\conv}{conv}

\DeclareMathOperator{\diam}{diam}
\DeclareMathOperator{\as}{as}
\DeclareMathOperator{\dell}{del}
\DeclareMathOperator{\divv}{div}
\DeclareMathOperator{\ldel}{ldel}
\DeclareMathOperator{\ldiv}{ldiv}

\newcommand{\dint}{\mathrm{d}}

\renewcommand{\phi}{\varphi}

\newtheorem {theorem}{Theorem}[section]

\newtheorem {example}[theorem]{Example}

\theoremstyle{definition}
\newtheorem {definition}{Definition}[section]
\newtheorem {remark}[theorem]{Remark}
\newtheorem {problem}[theorem]{Problem}

\title[Approximating Curved Shapes by Polytopes]{From Circles to Convex Bodies:\\ Approximating Curved Shapes by Polytopes}

\author{Steven Hoehner}
\date{\today}

\begin{document}

\setcounter{footnote}{0}
\maketitle

\begin{abstract}
Polytopes are the basic finite data structures for convex sets: they appear as feasible regions in
linear optimization, as geometric summaries in algorithms, and as random objects in stochastic
geometry. A natural geometric question is therefore: how well can a smooth, curved convex body be
approximated by a polytope with only $N$ faces? A striking phenomenon is that in $\R^d$, many seemingly different approximation errors--such as volume, surface area, and others)
often decay like $N^{-2/(d-1)}$ when the body has smooth, positively curved boundary. This survey article
offers a guided tour of that ``universal exponent'', starting from the classical approximation of a
circle by an $N$-gon and building intuition via spherical caps and curvature. We then survey a few
representative theorems--including results showing that random polytopes can be almost as good as
best possible ones--and explain why the Euclidean ball is a natural benchmark for the ``hardest" case. We also highlight a recently introduced projection-based distance that compares bodies through the average
distance of their shadows. Finally, we list accessible open problems about sharp constants, dimension
dependence, and gaps between known upper and lower bounds.
\end{abstract}

\renewcommand{\thefootnote}{}
\renewcommand{\thefootnote}{\arabic{footnote}}
\setcounter{footnote}{0}

\tableofcontents

\section{A warm-up: polygons and the circle}\label{sec:warmup}

Approximating a circle by a regular $N$-gon is a classical problem in geometry dating back to the ancients. If we need to estimate the area $\vol_2(B_2^2)=\pi$ of the unit disk $B_2^2$ by a regular $N$-gon $P_N$ inscribed in its boundary $\partial B_2^2=\mathbb{S}^1$, how close can we get for a fixed $N$? 
A first observation is that the relevant errors decay \emph{quadratically in $1/N$}.

\begin{example}\label{ex:N2}
Let $P_N$ be a regular  $N$-gon inscribed in the unit disk.
\begin{enumerate}[(i)]
\item The Hausdorff distance satisfies $d_H(B_2^2,P_N)\asymp N^{-2}$.
Indeed, the maximum radial deficit occurs halfway between consecutive vertices, where
$1-\cos(\pi/N)\sim \frac{\pi^2}{2N^2}$.
\item The area deficit satisfies $\vol_2(B_2^2)-\vol_2(P_N)\asymp N^{-2}$.
A standard computation gives $\vol_2(P_N)=\frac{N}{2}\sin(2\pi/N)$ and
$\pi-\frac{N}{2}\sin(2\pi/N)\sim \frac{2\pi^3}{3N^2}$.
\end{enumerate}
\end{example}

\begin{remark}[Key takeaway]
    In the plane, ``flat pieces'' approximate a curved boundary with a \emph{quadratic error} in the
typical mesh size. 
\end{remark}

As we will see next, the higher-dimensional exponent $2/(d-1)$ comes from combining this quadratic
behavior with the fact that $N$ faces divide a $(d-1)$-dimensional boundary into patches of typical
diameter about $N^{-1/(d-1)}$.

\vspace{2mm}

\begin{center}
    \includegraphics[scale=0.7]{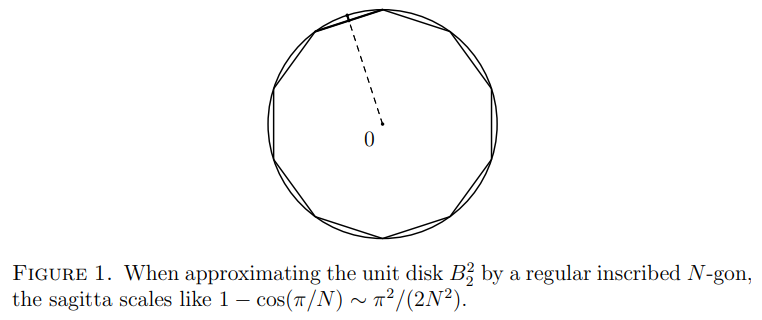}

\end{center}

\section{Why the exponent $2/(d-1)$ appears}\label{sec:exponent}

Let $K\subset\R^d$ be a convex body, i.e., a convex, compact set with nonempty interior. For a fixed integer $N\geq d+1$, let $\mathcal{Q}_N$ denote a family of polytopes--possibly depending on $K$--of ``complexity" at most $N$  (e.g., at most $N$ vertices, facets, intermediate faces, or flags). In algorithms, replacing $K$ by $P$ turns continuous optimization or geometry into finite linear constraints. It is thus  desirable to quantify the optimal error that can be achieved at a given level of complexity. The classical \emph{best approximation problem} asks for the asymptotic behavior, as $N\to\infty$, of the quantity
\begin{equation}\label{eq:inf}
\inf_{P\in\mathcal{Q}_N}d(K,P),
\end{equation}
where $d(K,P)$ is some measure of the ``distance" between $K$ and $P$. The probabilistic analogue of \eqref{eq:inf} replaces the infimum by the typical size (e.g., expected value) of $d(K,P_N)$ for a suitable random polytope model $P_N$. 

Starting from Problem \eqref{eq:inf}, the subject of polytopal approximation of convex bodies has three intertwined axes:
\begin{itemize}
    \item {\bf Metric axis}: Various ``distances" can and have been considered, including, for example: the Hausdorff metric, symmetric difference metric, surface area deviation, and intrinsic volume metric.

    \item {\bf Positioning axis}: The relative position of the body and the polytope may be, e.g., \emph{inscribed} ($P\subset K$), \emph{circumscribed} ($P\supset K$) or \emph{arbitrary position} (no containment).

    \item {\bf Complexity axis}: The approximating polytopes may have, e.g., a bounded number of vertices, facets, $k$-dimensional faces, or flags. 
\end{itemize}

Along these axes, a recurring theme in the area is that for a $C_+^2$ convex body $K\subset\R^d$ (that is, $K$ has boundary $\partial K$ of differentiable class at least $C^2$ and strictly positive Gauss--Kronecker curvature $\kappa_K$), many natural distances exhibit the same first-order decay, namely,
\begin{equation}\label{main-eqn}
    \inf_{P\in \mathcal{Q}_N}d(K,P)\asymp\NN,\qquad \E[d(K,P_N)]\asymp \NN
\end{equation}
as $N\to\infty$. The leading constants in such asymptotic formulas are governed by curvature integrals--including the affine surface area and its variants--depending only on $K$. This is classical for the symmetric difference metric in the inscribed and circumscribed settings (say), and the general asymptotic formula \eqref{main-eqn} will be the reference point for everything that follows. This section explains, with a heuristic computation, why $N^{-2/(d-1)}$ appears so often as the optimal rate in the polytopal approximation of smooth, strictly convex bodies in $\R^d$.

\subsection{A semi-rigorous cap computation on the sphere}

For a simple explanation showing how the exponent $2/(d-1)$ is attained, we will select the Euclidean unit ball $B_2^d=\{(x_1,\ldots,x_d)\in\R^d: x_1^2+\ldots+x_d^2\leq 1\}$.
Consider approximating the ball's boundary $\partial B_2^d$ near a point by a supporting hyperplane.
A standard model for the local geometry of $\partial B_2^d$ is a spherical cap centered at the supporting hyperplane's point of tangency; without loss of generality, assume that this point is $e_d=(0,\ldots,0,1)$. Fix $h\in(0,1)$ and consider the spherical cap of height $h$ defined by
\[
C(h):=\{x\in B_2^d:\ x_d\geq 1-h\}.
\]
Its base is the $(d-1)$-dimensional ball in the hyperplane $x_d=1-h$ with center $(1-h)e_d$ and radius
\[
r(h)=\sqrt{1-(1-h)^2}=\sqrt{2h-h^2}\sim \sqrt{2h}.
\]
Thus the cap sits above a patch of boundary whose typical radius is $\asymp h^{1/2}$, and hence whose 
typical diameter is $\asymp h^{1/2}$. The cap's volume is
\[
\vol_d(C(h))=\int_{1-h}^1 \vol_{d-1}\!\bigl(\sqrt{1-t^2}\cdot B_2^{d-1}\bigr)\,\dint t
= \vol_{d-1}(B_2^{d-1})\int_{1-h}^1 (1-t^2)^{\frac{d-1}{2}}\,\dint t.
\]
Substitute $t=1-s$ with $s\in[0,h]$. Then $1-t^2=1-(1-s)^2=2s-s^2\sim 2s$ for small $s$, so
\[
\vol_d(C(h))
\sim \vol_{d-1}(B_2^{d-1})\int_0^h (2s)^{\frac{d-1}{2}}\,\dint s
= \vol_{d-1}(B_2^{d-1})\cdot\frac{(2h)^{\frac{d+1}{2}}}{d+1}.
\]
In short,
\begin{equation}\label{eq:cap-volume-scaling}
\vol_d(C(h))\asymp h^{\frac{d+1}{2}}.
\end{equation}

To pass from caps to the asymptotic order $N^{-2/(d-1)}$, note that
if a polytope approximates the ball by ``supporting patches'' of typical boundary diameter $\rho$,
then the associated cap height scales like $h\asymp \rho^2$ (a paraboloid law), and the associated
volume error per patch scales like $h^{(d+1)/2}\asymp (\rho^2)^{(d+1)/2}=\rho^{d+1}$.
If $N$ patches cover the boundary, then a typical patch has $(d-1)$-dimensional area $\asymp 1/N$ and
hence diameter $\rho\asymp N^{-1/(d-1)}$. Therefore, a typical per-patch error is
$\rho^{d+1}\asymp N^{-(d+1)/(d-1)}$, and multiplying by $N$ patches gives a total error of order
$N\cdot N^{-(d+1)/(d-1)}=N^{-2/(d-1)}$.

\begin{remark}[Key takeaway]
The universal exponent $2/(d-1)$ is the product of two facts:
(i) smooth curvature produces a quadratic deviation from a tangent hyperplane, and
(ii) $N$ faces distribute approximation effort across a $(d-1)$-dimensional boundary.
\end{remark}

\vspace{2mm}

\begin{center}
    \includegraphics[scale=0.7]{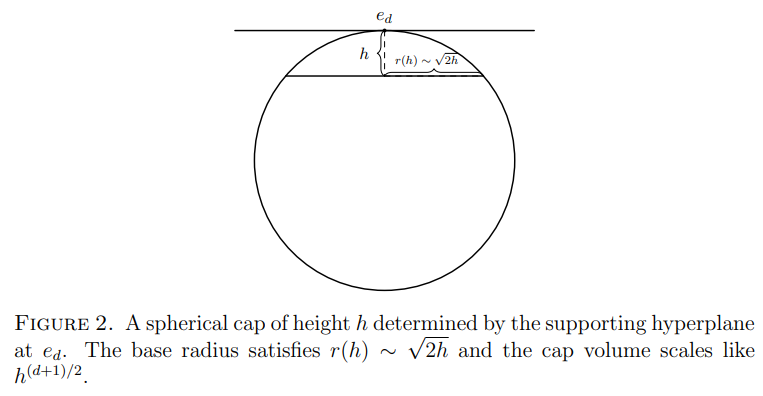}
\end{center}

\section{A short timeline of how the subject developed}\label{sec:timeline}

 In this section, we present a brief  roadmap of how the main ideas entered the subject.  This is by no means a complete history, but it highlights several milestones that shaped the modern theory. 

\begin{itemize}
\item \textbf{1950s--1960s:} Early work began with the approximation of bodies in $\R^2$ and $\R^3$ by L. Fejes T\'oth \cite{FejesToth}. The approximation properties of ellipsoids and Euclidean balls were also studied, leading to ``extremal body'' principles (e.g.,  \cite{Macbeath}).

\item \textbf{1963:} R\'enyi and Sulanke \cite{RenyiSulanke} initiated the systematic study of random polytopes in the plane. They analyzed the expected missing area of a random inscribed polygon, revealing how the expectation depends on the curvature of $K$.

\item \textbf{1974--1975:} Dudley \cite{Dudley1974,Dudley1979} and  Bronshtein and Ivanov \cite{Bronshtein-Ivanov} proved uniform bounds: every convex body can be approximated to Hausdorff accuracy $\varepsilon$ by a polytope with a controlled complexity: Dudley's bound requires only $O((\diam(K)/\varepsilon)^{(d-1)/2})$ facets, whereas Bronshteyn and Ivanov's bound requires the same number of vertices.

\item \textbf{1980s--2000s:} B\'ar\'any, B\"or\"oczky, Gruber, Ludwig, M\"uller, Sch\"utt, Werner  and many others developed the high-dimensional best and random approximation theories. Their results clarified how curvature functionals, such as the affine surface area, govern the sharp asymptotics for the error in best approximation. For random polytopes, their results show how the optimal sampling density is curvature-dependent, and that the resulting
expected symmetric difference error has the optimal $N^{-2/(d-1)}$ decay rate. For some examples, we refer the reader to \cite{Affentranger1991,Barany1992,Barany-Larman,glasgrub,GRS97,GRS98,Gruber88,GruberOut,Ludwig1999,muller1990,SW2003}.

\item \textbf{2000s--present:} In the past couple of decades, the field has seen renewed emphasis on sharp constants and optimal dimension dependence, understanding approximation without containment restrictions, and on extending the theory beyond vertices and facets to intermediate face numbers and flag constraints. Examples include: sharp surface area deviation results for the ball and smooth bodies \cite{BoroczkyCsikos,HSW,Kur2017,GTW2021}; asymptotic estimates while bounding intermediate face numbers  \cite{Boroczky2004,Boroczky-Fodor-Vigh,Boroczky-Gomis-Tick,HSW-2025}; concentration inequalities and Central Limit Theorems for random polytopes \cite{HK-DCG,reitzner2003random,vu2005sharp,vu2006central}; and sharp  approximation of the ball by polytopes with an exponential number of facets \cite{Kur2017,LSW}.
\end{itemize}

\smallskip

Other surveys on the subject can be found in, e.g.,  \cite{BaranySurvey,HugSurvey,PSW-2022,ReitznerSurvey} and \cite[Chapter 11]{GruberBook}.

\begin{remark}[Key takeaway]
The modern theory combines \emph{geometry} (curvature and caps), \emph{combinatorics} (how faces
partition the boundary), and \emph{probability} (random polytopes as almost-optimal constructions).
\end{remark}

\section{Two classical metrics: Hausdorff and symmetric difference}\label{sec:classical-metrics}

There are many reasonable ways to quantify how close $K$ and $P$ are. Two popular choices are:
\begin{itemize}
\item the \emph{Hausdorff metric} \[d_H(K,P)=\max\left\{\sup_{x\in K}d(x,P),\sup_{y\in P}d(y,K)\right\},\] where $d(x,K)=\inf_{y\in K}\|x-y\|_2$, and 

\item the \emph{symmetric difference metric} 
\[
d_S(K,P)=\vol_d(K\triangle P)=\vol_d(K)+\vol_d(P)-2\vol_d(K\cap P), 
\]
where $K\triangle P=(K\cup P)\setminus(K\cap P)$. Note that when $P\subset K$ (respectively, $P\supset K$), $d_S(K,P)$ reduces to the \emph{volume difference} $\vol_d(K)-\vol_d(P)$ (respectively, $\vol_d(P)-\vol_d(K)$).
\end{itemize}
These distances emphasize different phenomena: the Hausdorff distance is sensitive to the single worst direction, while
the symmetric difference metric measures  distance by averaging error across the whole boundary.

\subsection{Hausdorff approximation}\label{subsec:dudley}

A cornerstone result, proved independently by Dudley \cite{Dudley1974,Dudley1979} and by Bronshtein and Ivanov \cite{Bronshtein-Ivanov}, says that \emph{any} convex body admits a polyhedral approximation with $N$ facets (or $N$ vertices) that has a Hausdorff error of about $N^{-2/(d-1)}$. This is the same exponent as in Section~\ref{sec:exponent}, but now it arises from a different mechanism, namely, metric entropy.

\begin{theorem}[Hausdorff approximation by polytopes]\label{thm:dudley}
Let $K\subset\R^d$ be any convex body.
Then for every $\varepsilon\in(0,1]$:
\begin{itemize}
    \item[(a)]  There exists a polytope $P$ in $\R^d$ with at most 
$(\diam(K)/\varepsilon)^{(d-1)/2}$ facets such that
\[
K\subset P\subset K+\varepsilon B_2^d.
\]
Equivalently, $d_H(K,P)\leq \varepsilon$.

\item[(b)] There exists a polytope $Q$ in $\R^d$ with at most 
$(\diam(K)/\varepsilon)^{(d-1)/2}$ vertices such that
\[
K-\varepsilon B_2^d\subset Q\subset K.
\]
Equivalently, $d_H(K,Q)\leq \varepsilon$.
\end{itemize}
In particular, there exists a constant $C(K)>0$ such that for every $N\geq 1$, we can find a polytope
$P$ with at most $N$ facets (or at most $N$ vertices) satisfying
\[
d_H(K,P)\leq C(K)\NN.
\]
Moreover, the exponent $(d-1)/2$ is best possible in the worst case (for example, when $K=B_2^d$), up to the value of the constant.
\end{theorem}

For a modern treatment of this result and ones similar to it, see \cite{AryaDaFonsecaMount2012Mahler,har-peled}. This result standard in computational geometry, where it is used in the construction of geometric algorithms (see, e.g., \cite{AryaDaFonsecaMount2012Mahler}). Curvature enters only in refinements of the constant; sharp curvature-driven asymptotics for inscribed or circumscribed polytopes require smoothness. The precise asymptotics are governed by the rate $\NN$. 

The next result, due to Gruber \cite[Theorem 3]{GruberI}, describes the precise asymptotic behavior of the best Hausdorff approximation. In the sequel, for a convex body $K\subset\R^d$ we let $\mathcal{P}_N^{\rm in}(K)$ (respectively, $\mathcal{P}_N^{\rm out}(K)$) be the set of all polytopes $P\subset K$ ($P\supset K$) such that $P$ has at most $N$ vertices (facets). Also, let $\kappa_K(x)$ denote the Gauss--Kronecker curvature at $x\in\partial K$, and let $\mu_{\partial K}$ be the surface area measure on $\partial K$.

\begin{theorem}[Asymptotic Hausdorff approximation]\label{thm:inscribed-hausdorff}
Let $K\subset\R^d$ be a convex body of class $C^2_+$, and let $\mathcal{Q}_N$ be either $\mathcal{P}_N^{\rm in}(K)$ or $\mathcal{P}_N^{\rm out}(K)$. Then the best Hausdorff approximation  satisfies the sharp asymptotic formula
\[
\inf_{P\in\mathcal{Q}_N} d_H(K,P)
\sim
\frac12
\left(
\frac{\vartheta_{d-1}}{\vol_{d-1}(B_{d-1})}
\int_{\partial K} \kappa_K(x)^{1/2}\, \dint\mu_{\partial K}(x)
\right)^{\frac{2}{d-1}}
\NN
\qquad \text{as }N\to\infty,
\]
where $\vartheta_{d-1}$ is the minimal
covering density of $\R^{d-1}$ by unit balls. In particular, for $K=B_2^d$ we have  $\inf_{P\in\mathcal{Q}_N}d_H(B_2^d,P)
\sim c_d\NN$ as $N\to\infty$, with an explicit dimensional constant $c_d>0$.
\end{theorem}

\begin{remark}
The curvature weight $\kappa_K^{1/2}$ in Theorem~\ref{thm:inscribed-hausdorff} says, informally, that an optimal   polytope should place its faces (i.e., vertices or facets) more densely where the boundary is more curved: in the large $N$ limit, the vertices become asymptotically equidistributed along $\partial K$ with density proportional to $\kappa_K^{1/2}$, which means regions of higher curvature receive a larger share of the available faces.
\end{remark}

\begin{remark}[Key takeaway]
Hausdorff approximation is governed by ``worst" directions. Even for non-smooth bodies one gets the same uniform rate, while smoothness is needed to state sharp curvature-dependent limits.
\end{remark}

\subsection{Symmetric difference and curvature}\label{subsec:symmdiff}

For the symmetric difference metric, smooth curvature enters the picture in a more direct way. When $K$ is strictly convex with $C_+^2$ boundary, the optimal approximation rate still behaves like $N^{-2/(d-1)}$, but the leading constant depends on the curvature of $K$ through its affine surface area.

\begin{definition}[Affine surface area]\label{def:as}
Let $K$ be a convex body in $\R^d$. The  \emph{affine surface area} of $K$ is defined by
\[
\as(K)=\int_{\partial K}\kappa_K(x)^{\frac{1}{d+1}}\,\dint \mu_{\partial K}(x).
\]
\end{definition}
\noindent You can think of $\as(K)$ as a ``curvature-weighted" notion of surface area, where regions of higher curvature contribute more to $\as(K)$.

\vspace{2mm}

In the classical inscribed ($P\subset K$) and circumscribed ($P\supset K$) models, precise limit formulas are known. 
The following formula describes how curvature appears in the first-order sharp constants for best approximation.

\begin{theorem}[Prototype asymptotic formula]\label{thm:prototype}
Let $K\subset\R^d$ be a  $C^2_+$ convex body. In several standard polytope models $\mathcal{Q}_N$,  
the optimal approximation satisfies the following asymptotic formula:
\[
\inf_{P\in\mathcal{Q}_N} \vol_d(K\triangle P) \asymp C(K)\NN.
\]
\end{theorem}

\noindent Here $C(K)$ is a weighted curvature integral (e.g., $C(K)=c_d\int_{\partial K}\kappa_K(x)^\alpha\,\dint\mu_{\partial K}(x)$)  whose weight (exponent $\alpha$) depends on $\mathcal{Q}_N$ and $d$. In fact, in each model the ratio tends to an explicit constant; see \eqref{eq:inscribed-best} and Remark \ref{rmk:circumscribed} below.  Moreover, the sharp limits contain dimension-dependent constants connected to optimal local tilings, triangulations or power diagrams.

Let us illustrate these principles through a few examples. Gruber \cite{Gruber88} proved that if $K\subset\R^d$ is a $C_+^2$ convex body, then 
\begin{equation}\label{eq:inscribed-best}
\inf_{P\in\mathcal{P}_N^{\rm in}(K)}\left(\vol_d(K)-\vol_d(P)\right)\sim \frac{\dell_{d-1}}{2}\cdot\as(K)^{\frac{d+1}{d-1}}\NN\qquad \text{as }N\to\infty.
\end{equation}
Here $\dell_{d-1}$ is the \emph{Delone triangulation number}, a positive constant whose value depends only on the dimension. This number is connected with Delone triangulations in $\R^{d-1}$, which arise in the optimal local approximation of boundary patches of $K$. The exact value of $\dell_{d-1}$ is known only for dimensions 2 and 3. Fejes T\'oth \cite{FejesToth} showed that $\dell_1=1/6$, and Gruber \cite{Gruber88} proved that $\dell_2=1/(2\sqrt{3})$. Mankiewicz and Sch\"utt \cite{MaS1,MaS2} proved that $\dell_{d-1}=\frac{d}{2\pi e}\left(1+O\left(\tfrac{\ln d}{d}\right)\right)$ as $d\to\infty$.

\begin{remark}\label{rmk:circumscribed}
    For circumscribed polytopes with at most $N$ facets, Gruber \cite{GruberOut} proved that the optimal approximation is again given by \eqref{eq:inscribed-best}, but with $\dell_{d-1}$ replaced by $\divv_{d-1}$, the \emph{Delone--Voronoi tiling number}. These numbers depend only on the dimension $d$ and are connected with Voronoi tilings in $\R^{d-1}$. Only the values $\divv_1=1/12$ (see \cite{FejesToth}) and $\divv_2=5/(18\sqrt{3})$ (see \cite{GruberOut}) are known explicitly. It is known that $\divv_{d-1}=\frac{d}{2\pi e}\left(1+\tfrac{\ln d}{d}+O\left(\tfrac{1}{d}\right)\right)$ as $d\to\infty$  \cite{HK-DCG,HSW-2025}.

\end{remark}

\subsubsection{Polytopes in arbitrary position}

If we drop the containment restriction $P\subset K$ or $P\supset K$, how much can we improve the estimate? This question is especially interesting because containment forces every face of $P$ to lie on one side of $\partial K$, so the error is essentially one-sided, i.e., you only ``miss'' volume or
only ``add'' volume. In arbitrary position, however, the polytope may locally overshoot and undershoot
$\partial K$ in different regions, allowing a more efficient global trade-off. As a consequence,
while the universal rate $N^{-2/(d-1)}$  persists for $C_+^2$ bodies, by dropping the containment restriction the sharp constants can improve
substantially---in some cases by a factor of dimension---and new phenomena arise that have no analogue in the inscribed or circumscribed settings.

In the groundbreaking work \cite{Ludwig1999}, Ludwig derived analogues of the asymptotic formula \eqref{eq:inscribed-best} and its circumscribed counterpart. To state the result, let $K\subset\R^d$ be a $C_+^2$ convex body, and let $\mathcal{P}_{d,N}$ denote the set of all polytopes in $\R^d$ with at most $N$ vertices. Ludwig \cite{Ludwig1999} proved that 
    \begin{equation}\label{eq:arbitrary-best}
    \inf_{P\in\mathcal{P}_{d,N}}\vol_d(K\triangle P)\sim\frac{\ldel_{d-1}}{2}\cdot\as(K)^{\frac{d+1}{d-1}}\NN\qquad \text{as }N\to\infty,
    \end{equation}
where  $\ldel_{d-1}$ is a constant depending only on $d$ called the  \emph{Laguerre--Delone triangulation number}. It is known that $\ldel_1=1/16$ (see \cite{Ludwig1999}) and $\ldel_2=1/(6\sqrt{3})-1/(8\pi)$ (see \cite{Boroczky-Ludwig}). In general dimensions, it is known that
\begin{equation}\label{ldel-estimates}
    \frac{c}{d} \leq \ldel_{d-1} \leq C
\end{equation}
for some positive absolute constants $c,C$. The lower estimate is due to B\"or\"oczky \cite{Boroczky2004}, while the upper estimate is due to Ludwig, Sch\"utt and Werner \cite{LSW}. These estimates show that dropping the containment restriction improves the estimate by at least a factor of dimension. This happens because when you allow overshoot and undershoot, cancellations become possible, which one-sided approximation forbids. Closing the gap between the upper and lower estimates for $\ldel_{d-1}$ is one of the major open problems in the area. 

\begin{remark}
    Ludwig \cite{Ludwig1999} also established a similar asymptotic formula for the best approximation of $K$ by arbitrarily positioned polytopes with at most $N$ facets, with dimensional constant $\ldiv_{d-1}$, the \emph{Laguerre--Dirichlet--Voronoi tiling number}. For these numbers, it is known that $\ldiv_1=1/16$ (see \cite{Ludwig1999}) and $\ldiv_2=5/(18\sqrt{3})-1/(4\pi)$ (see \cite{Boroczky-Ludwig}). In high dimensions, sharp estimates have been obtained:
\begin{equation}\label{ldiv-eqn}
c_1 \leq \ldiv_{d-1} \leq c_2.
\end{equation}
The lower estimate is due to Ludwig, Sch\"utt and Werner \cite{LSW}, and the upper estimate is due to Kur \cite{Kur2017}, who also derived numerical estimates for $c_1,c_2$. 
\end{remark}

It is natural to ask why there is a tighter estimate for $\ldiv_{d-1}$ than $\ldel_{d-1}$. In the general position case, if $P$ is a best-approximating polytope with a restricted number of facets, then exactly half of the $(d-1)$-volume of each facet of $P$ lies inside $K$ \cite{LSW}. This condition translates to strong balancing of the facets in the overshoot/undershoot trade-off.

\begin{remark}[Key takeaway]
For smooth strictly convex bodies, volume-based approximation is governed by two universal features:
(i) the error decays like $N^{-2/(d-1)}$, reflecting the local quadratic (paraboloid) nature of
a smooth boundary, and (ii) the leading constant splits into a \emph{geometric part} and a \emph{combinatorial part}.  The geometric part is a weighted  curvature integral, which tells us where an optimal polytope should ``spend'' its faces, while
the combinatorial part is a dimensional constant coming from optimal local tilings, triangulations or power diagrams 
(i.e., $\dell_{d-1},\divv_{d-1}$ in the inscribed/circumscribed settings and
$\ldel_{d-1},\ldiv_{d-1}$ in arbitrary position).  Dropping the containment restriction does not change the
exponent $2/(d-1)$, but it can improve the sharp constant substantially; in particular, for facet   
approximation, the best-known bounds for $\ldiv_{d-1}$ imply a dimension-level advantage.
\end{remark}

\section{Random polytopes can be (nearly) as good as best possible}\label{sec:random}

One of the most remarkable aspects of this subject is that, in an asymptotic sense, probabilistic constructions often produce
polytopes whose average error matches the best possible order in $N$. Let us illustrate this phenomenon with an example. 

\subsection{Convex hull of random boundary points}

Given a convex body $K\subset\R^d$, sample  $N$ random points $X_1,\dots,X_N$  from  $\partial K$ independently  according to a continuous, positive probability density $f$, and let $P_f=\conv\{X_1,\dots,X_N\}$ denote the convex hull formed by these points. Under some regularity conditions on $K$, Sch\"utt and Werner \cite{SW2003} proved that
\begin{equation}\label{SW-random}
   \lim_{N\to\infty}\NNN\E[d_S(K,P_f)]= \lim_{N\to\infty}\NNN\left(\vol_d(K)-\E[\vol_d(P_f)]\right)=c_d\int_{\partial K}\frac{\kappa_K(x)^{\frac{1}{d-1}}}{f(x)^{\frac{2}{d-1}}}\,\dint\mu_{\partial K}(x),
\end{equation}
where $c_d$ is an explicit positive constant, and that the minimizing density is 
\begin{equation}\label{min-density}
f_{\min}(x)=\frac{\kappa(x)^{\frac{1}{d+1}}}{\as(K)},\quad x\in\partial K.
\end{equation}
Choosing the optimal density $f_{\min}$, the right-hand side of \eqref{SW-random} becomes $c_d\as(K)^{\frac{d+1}{d-1}}$. From this we get 
\begin{equation}
    \lim_{N\to\infty}\frac{\E[d_S(K,P_{f_{\min}})]}{\displaystyle\inf_{P\in\mathcal{P}_N^{\rm in}(K)}d_S(K,P)}=\frac{c_d}{\dell_{d-1}/2}\sim 1\qquad\text{as }d\to\infty,
\end{equation}
where we used the estimate $\frac{1}{2}\dell_{d-1}\leq c_d\leq \frac{1}{2}\dell_{d-1}(1+O(\frac{\ln d}{d}))$ (see \cite{SW2003}).


More recently, random polytopes have been used to study the arbitrary position case. As in the inscribed and circumscribed cases, it turns out that dropping the containment restriction again improves the estimate by a factor of dimension. 

\begin{theorem}[Grote--Werner \cite{GW2018}]\label{thm:GW}
Let $K\subset\R^d$ be a $C^2_+$ convex body, and let $P_f$ be the convex hull of $N$ random  points sampled from $\partial K$ independently according to a positive, continuous probability density $f$. Then for all sufficiently large $N$,
\[
\E[\vol_d(K\triangle P_f)] \leq C\NN\int_{\partial K}\frac{\kappa_K(x)^{\frac{1}{d-1}}}{f(x)^{\frac{2}{d-1}}}\,\dint\mu_{\partial K}(x)
\]
where $C>0$ is an absolute constant. The minimizing density is given by $f_{\min}$ defined in \eqref{min-density}.
\end{theorem}


\begin{remark}[Key takeaway]
The random polytope model does more than just provide computational tools from probability, or yield a convenient upper bound: it reveals an
explicit variational principle for how vertices should be distributed.
Formula~\eqref{SW-random} shows that the expected missing volume is governed by a curvature density
trade-off, and the minimizing density $f_{\min}$ tells us that an asymptotically optimal random polytope samples more frequently where $\partial K$ is more curved. A parallel picture holds for random circumscribed polytopes built as intersections of random supporting halfspaces: with an appropriate choice of directional distribution, one again obtains the universal rate $N^{-2/(d-1)}$, and the resulting random construction is asymptotically near-optimal compared to the best possible circumscribed polytope; see, for example, \cite{BHK} and the references therein.
\end{remark}


\subsection{Surface area approximation: a closely related story}

Along the ``metric axis", one can measure the error in approximation with respect to other types of volume. The \emph{surface area deviation}, defined by
\[
\Delta_s(K,P):=\vol_{d-1}(\partial K)+\vol_{d-1}(\partial P)-2\vol_{d-1}(\partial(K\cap P)),
\]
behaves much like the symmetric volume difference, but is more sensitive to boundary geometry, rather than interior
volume. Although $\Delta_s$ is not a metric--it fails the triangle inequality; see \cite{BHK}--it is still a natural and useful notion of error. It plays the same role for surface area that the symmetric difference metric plays for volume, in that it measures how much ``boundary mass'' is lost in $K\cap P$ compared to the two bodies separately. Moreover, in polytopal approximation it
admits sharp upper and lower bounds that exhibit the same universal rate $N^{-2/(d-1)}$ for smooth
strictly convex bodies, but with different sharp constants than the volume-based theory. For some
examples, we refer the reader to  \cite{BoroczkyCsikos,GTW2021,HSW,Kur2017}.

For the Euclidean ball, sharp upper and lower bounds for $\Delta_s(B_2^d,P)$ for generally
positioned polytopes were established in \cite{HSW,Kur2017}. For a $C_+^2$ convex body $K\subset\R^d$, Grote, Th{\"a}le and Werner \cite{GTW2021} proved optimal-order upper bounds--again, of order
$N^{-2/(d-1)}$--with leading terms connected to the \emph{$d$-affine surface area}
\[
\as_d(K)=\int_{\partial K}\frac{\kappa_K(x)^{1/2}}{\langle x,N_K(x)\rangle^{\frac{d-1}{2}}}\,\dint\mu_{\partial K}(x),
\]
where $N_K(x)$ is the unique outer unit normal of $K$ at $x$.

\begin{theorem}[Grote--Th{\"a}le--Werner \cite{GTW2021}]\label{thm:GTW}
Let $K\subset\R^d$ be a $C^2_+$ convex body. There exists an integer $N_K$ depending only on $K$ such that for all $N\geq N_K$, there exists a polytope $P$ in $\R^d$ with at most $N$ vertices such that
\[
\Delta_s(K,P)\leq Cd\as_d(K)^{\frac{2}{d-1}}\vol_{d-1}(\partial K)\NN
\]
where $C$ is a positive absolute constant.
\end{theorem}

\subsubsection{The mean width deviation}

Another way to measure ``volume-type'' error is via \emph{mean width}, which is built from the support
function and is therefore closer in spirit to the Hausdorff metric than to volume difference.
The \emph{mean width} of a convex body $K$ is
\[
w(K)=\frac{2}{\vol_{d-1}(\partial B_2^d)}\int_{\partial B_2^d}h_K(u)\,\dint\mu_{B_2^d}(u),
\qquad
h_K(u)=\sup_{x\in K}\langle x,u\rangle .
\]
Geometrically, $h_K(u)$ is the distance from the origin to the supporting hyperplane of $K$ with outer
normal $u$, and $w(K)$ is the average of these support values over all directions.  In the plane one
has the pleasing identity $w(K)=\pi^{-1}\vol_1(\partial K)$, so mean width can be viewed as a
higher-dimensional analogue of perimeter.

To compare $K$ with a polytope $P$, one may use the \emph{mean width deviation}
\[
\Delta_w(K,P):=w(K)+w(P)-2w(K\cap P).
\]
Although $\Delta_w$ is not a metric, it inherits the same ``intersection'' structure as the volume and surface area deviations, and is a natural error functional along the metric axis: it detects how much
``directional size'' is lost in the intersection $K\cap P$.  Moreover, mean width fits neatly into
the intrinsic volume hierarchy discussed later in Section~\ref{sec:shadows}: $V_1(K)$ is a constant multiple of $w(K)$, so $\Delta_w$ is (up to a
fixed dimensional factor) the $j=1$ case of the intrinsic volume deviation
$V_j(K)+V_j(P)-2V_j(K\cap P)$.

Asymptotically, when $K$ is $C^2_+$, the best (inscribed or circumscribed) mean width approximation
again exhibits the universal rate
\[
\inf_{P\in\mathcal{Q}_N}\Delta_w(K,P)\sim c_w(K)\,\NN\qquad\text{as }N\to\infty,
\]
where the leading constant $c_w(K)$ is given by a curvature integral, but with a \emph{different}
curvature weight than in the volume or surface area cases.  Consequently, the ``optimal distribution'' of vertices
(or facet normals) for mean width approximation is generally different from the optimal distribution
for volume and surface area approximations: the geometry (curvature) is the same, but the metric determines how that curvature is weighted.  For examples and sharp asymptotic formulas we refer to
\cite{Muller1989,glasgrub,Ludwig1999,BHK} and the references therein.

\begin{remark}[Key takeaway]
Mean width provides a bridge between Hausdorff-type information (support functions in each
direction) and volume-type information (since $w$ is proportional to $V_1$).  The same universal
rate $\NN$ appears again for smooth bodies, but the curvature weight--and hence the
asymptotically optimal ``face density''--depends on whether the error is quantified by volume, surface
area, or mean width.
\end{remark}

\section{Why the ball is a benchmark for the ``hardest" case}\label{sec:ball}

In many approximation problems, an optimal polytope spends more faces where curvature is high.
If curvature varies across $\partial K$, this provides flexibility. The Euclidean ball has constant
curvature, leaving no advantageous region, i.e., approximation effort must be spread uniformly.
This intuition is reflected in several extremal statements in the literature, beginning with a
classical theorem of Macbeath \cite{Macbeath}. It states that, among all convex bodies of a given volume,  ellipsoids are hardest to approximate by inscribed polytopes with a fixed number of vertices.     

\begin{theorem}[Macbeath's theorem]
For every convex body $K\subset\R^d$,
\begin{equation}\label{vol-worst-case}
    \inf_{P\in\mathcal{P}_N^{\rm in}(K)}\frac{d_S(K,P)}{\vol_d(K)}\leq \inf_{P\in\mathcal{P}_N^{\rm in}(B_2^d)}\frac{d_S(B_2^d,P)}{\vol_d(B_2^d)}.
\end{equation}
\end{theorem}
The stated result for ellipsoids now follows by applying an affine transformation. An analogous result also holds for the mean width deviation \cite{Florian-1989}.

\begin{remark}[Key takeaway]
The ball is a ``stress test'' for approximation methods: if we cannot concentrate faces where
curvature is high, then we resort to optimizing globally. This is essentially why the ball often reveals the sharp dimension dependence of constants. 
\end{remark}

If $d$ is large and we require $N$ to be exponential in $d$, then a polytope with $N$ facets can already approximate the ball quite  well. A representative result is due to Kur \cite{Kur2017}.

\begin{theorem}[Kur \cite{Kur2017}]\label{thm:Kur}
There exists an absolute constant $C>0$ such that for every $d\geq 2$ and $N\geq 10^d$, there is a
polytope $P$ in $\R^d$ with at most $N$ facets satisfying both
\[
\vol_d(B_2^d\triangle P)\leq C\vol_d(B_2^d)\NN
\quad\text{and}\quad
\Delta_s(B_2^d,P)\leq C\vol_{d-1}(\partial B_2^d)\NN.
\]
\end{theorem}
These bounds are optimal up to absolute constants in the relevant regime--the same lower bounds hold with different absolute constants \cite{HSW,LSW}.

\section{Beyond volume: comparing bodies by their shadows}\label{sec:shadows}

So far we've measured distance in $\R^d$ (volume error) or on the boundary (surface area deviation).
A different--and surprisingly natural--idea is to compare shapes through their \emph{shadows}. In $\R^3$, for example, two solids may look similar from most directions even if they differ significantly in some localized region. This motivates projection-based distances.

\subsection{A projection-based distance: The intrinsic volume metric}

Volume, surface area and mean width are examples of the \emph{intrinsic volumes}\footnote{Intrinsic volumes are also known as \emph{quermassintegrals}. The term comes from the German word ``Querma\ss", which means ``measure of a either a cross-section or a projection" \cite{SchneiderBook}, and the suffix ``-integral" refers to an  average or total value.} of $K$, denoted $V_1(K),\ldots,V_d(K)$. These quantities represent specialized average volumes of a convex body's projections. More specifically, the intrinsic volumes can be expressed via \emph{Kubota's formula}
\begin{equation}\label{eq:Kubota}
V_j(K)=c_{d,j}\int_{\Gr(d,j)}\vol_j(K|E)\,\dint\nu_{d,j}(E),\qquad j\in\{1,\ldots,d\},
\end{equation}
where $K|E$ denotes the orthogonal projection of $K$ onto the $j$-dimensional subspace $E$ of $\R^d$, $\nu_{d,j}$ is the Haar probability measure on the Grassmannian manifold $\Gr(d,j)$ of all $j$-dimensional subspaces of $\R^d$, and $c_{d,j}>0$ is an explicit constant. In particular, $V_d(K)=\vol_d(K)$ is the volume, $2V_{d-1}(K)=\vol_{d-1}(\partial K)$ is the surface area, and $V_1(K)$ is a constant multiple of the mean width $w(K)$. For more background on the intrinsic volumes of a convex body, see  \cite{SchneiderBook}.

Along the metric axis of the theory, the error in approximation can be measured more generally via the intrinsic volumes. For $j\in\{1,\dots,d\}$, the \emph{$j$th intrinsic volume metric} is defined by
\[
\delta_j(K,L)
:=c_{d,j}\int_{\Gr(d,j)} \vol_j\bigl((K|E)\triangle(L|E)\bigr)\,\dint \nu_{d,j}(E),
\]
where $c_{d,j}$ is the same normalization constant as in  \eqref{eq:Kubota}. When $j=d$, this is just the symmetric difference volume $\vol_d(K\triangle L)$.
When $j<d$, it measures the average volume discrepancy of the $j$-dimensional shadows.

\begin{center}
    \includegraphics[scale=0.7]{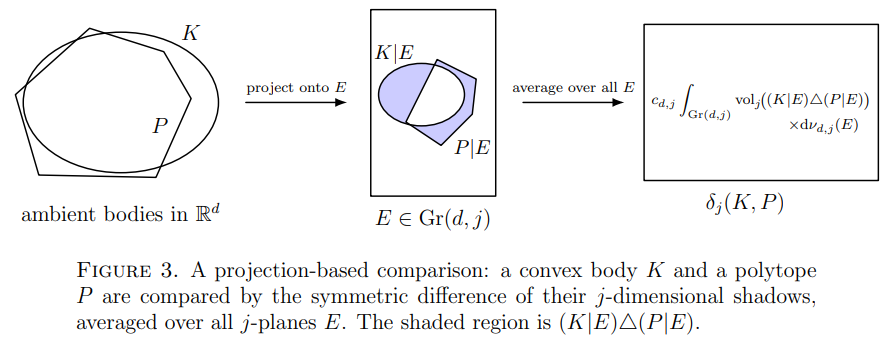}
\end{center}

The next result, which is from  \cite{BHK}, shows that the intrinsic volumes all provide essentially the same  geometric information in the approximation of the ball (see also \cite{BH-2024,HSW-2025}).

\begin{theorem}[Asymptotic intrinsic volume approximation]
    Let $\mathcal{Q}_N$ be either $\mathcal{P}_N^{\rm in}(B_2^d)$ or $\mathcal{P}_N^{\rm out}(B_2^d)$, and fix $j\in\{1,\ldots,d\}$. There exist positive absolute constants $C_1,C_2$ such that for all sufficiently large $N$,
    \[
C_1 jV_j(B_2^d)\NN\leq\inf_{P\in\mathcal{Q}_N}\delta_j(B_2^d,P)\leq C_2jV_j(B_2^d)\NN.
    \]
\end{theorem}
\noindent This estimate is asymptotically sharp in the sense that  the constants satisfy $c_1\sim c_2=\frac{1}{2}+O\left(\frac{\ln d}{d}\right)$ as $d\to\infty$ \cite{BHK}.

Remarkably, as shown in \cite{BHK}, (for each case of $\mathcal{Q}_N$) there exists a single polytope  which satisfies all $d$ of these inequalities, and thus is asymptotically optimal with respect to all intrinsic volumes  \emph{simultaneously}. In view of the extremal approximation properties of the ball, it is natural to ask: Is this asymptotic ``simultaneous approximation property" unique to the ball among all $C_+^2$ convex bodies?  Very recently, it was shown in \cite{Hoehner-2026} that the answer to this question is affirmative. The reason behind this phenomenon is simple: in order for a polytope to satisfy the volume ($j=d$) and mean width ($j=1$) bounds (say), the respective curvature-dependent vertex densities of the optimal polytopes must be equal, which forces the curvature to be constant.

For the intrinsic volume metric it is also known that, asymptotically, random inscribed polytopes achieve, on average, the same rate of approximation as the best polytope (see \cite{BH-2024,BHK}).


\begin{remark}[Key takeaway]
The intrinsic volume metrics $\delta_j$ measure \emph{typical} discrepancy rather than worst-case
discrepancy: instead of comparing $K$ and $P$ in the ambient space, they compare their
$j$-dimensional shadows and average over all directions. This makes $\delta_j$ robust to
localized differences and naturally links approximation to integral geometry through Kubota's
formula. For the Euclidean
ball, $\delta_j$ leads to a unified picture across all intrinsic volumes: the optimal error is
essentially of the same order for every $j$, and there are even polytopes that are asymptotically optimal for all $j$ simultaneously. For intermediate $j$ and for arbitrary position approximation, however, sharp lower bounds remain largely open.
\end{remark}

\section{A curated list of open problems}\label{sec:open}

We conclude with a short list of open problems phrased to be broadly accessible.

\begin{problem}[Sharp constants and dimension dependence]\label{prob:constants}
In several classical asymptotic estimates, the exponent $N^{-2/(d-1)}$ is known, but the sharp
dimension-dependent constants are mysterious in higher dimensions. Determine these constants, or
at least their precise growth as $d\to\infty$, and explain geometrically what structures
(e.g., triangulations or tilings) they encode.
\end{problem}

\begin{problem}[The $\ldel$ gap]\label{prob:ldel}
The best-known upper and lower estimates for $\ldel_{d-1}$ (see \eqref{ldel-estimates}) differ by a factor of dimension. Close this gap. Is there a dimension-free positive lower bound?
\end{problem}

\begin{problem}[Random vs. best for projection-based distances]\label{prob:delta}
Develop a sharp asymptotic approximation theory for the intrinsic volume metric $\delta_j$. For the Euclidean ball, determine the best possible decay of $\inf\{\delta_j(B_2^d,P): P\text{ has at most }N\text{ facets}\}$,
and determine whether the corresponding random polytope model achieves the same order.
\end{problem}

\begin{problem}[Intermediate face-number constraints]\label{prob:fk}
Most sharp results focus on polytopes with a fixed number of vertices  or facets. Extend the sharp asymptotic approximation theory to polytopes with a prescribed number of $k$-faces for $1\leq k\leq d-2$: what is the best possible accuracy, and what combinatorial types are extremal?
\end{problem}

\begin{problem}[Is the ball always the hardest?]\label{prob:hardest}
In several normalized settings, the Euclidean ball is the hardest body to
approximate for a given complexity. Determine for which metrics and models this remains true, especially in the circumscribed and arbitrary position models, and for surface area or shadow-based distances.
\end{problem}

\bibliographystyle{plain}
\bibliography{main}


\vspace{3mm}

\noindent {\sc Department of Mathematics \& Computer Science, Longwood University,\\ 201 High Street, Farmville, Virginia 23901}\\
\noindent {\it E-mail address:} {\tt hoehnersd@longwood.edu}

\end{document}